\documentclass[conference]{IEEEtran}
\usepackage{float}
\IEEEoverridecommandlockouts

\usepackage{pgfplots}
\usepackage{pgfplotstable}
\pgfplotsset{compat=newest}
\usepackage{cite}
\usepackage{amsmath,amssymb,amsfonts}
\usepackage{algorithm}
\usepackage{algpseudocode}  
\usepackage{float}  
\usepackage{graphicx}
\usepackage{textcomp}
\usepackage{xcolor}
\usepackage[font=footnotesize]{caption}
\usepackage{subcaption}
\usepackage{url}  
\usepackage{xurl}  
\usepackage{amsmath}
\usepackage{bm}  
\usepackage{tikz}
\usetikzlibrary{arrows.meta, positioning}

\def\BibTeX{{\rm B\kern-.05em{\sc i\kern-.025em b}\kern-.08em
    T\kern-.1667em\lower.7ex\hbox{E}\kern-.125emX}}
\begin{document}

\title{


Rapid Concurrent GPU--CPU Solvers for Scalable Unit Commitment in Large Power Grids

}
\author{\IEEEauthorblockN{Hussein Sharadga}
\IEEEauthorblockA{
\textit{~~~~~Texas A\&M International University}\\
Laredo, Texas, USA \\
hssharadga@tamu.edu}
\and
\IEEEauthorblockN{Yuhan Du, Javad Mohammadi}
\IEEEauthorblockA{
\textit{The University of Texas at Austin}\\
Austin, Texas, USA \\
{\{yuhandu, javadm\}@utexas.edu}}
}



\maketitle

\IEEEpubidadjcol

\begin{abstract}

This paper presents an accelerated solver for the unit commitment problem in large-scale power systems. The approach is based on the concurrent execution of GPU- and CPU-based optimization solvers on a single machine, with the solver that converges first terminating the other to minimize overall runtime. This strategy effectively harnesses the complementary strengths of different solvers. Convergence is further accelerated through a systematic and aggressive presolve approach. Numerical experiments on a 6,049-bus system with millions of decision variables and constraints demonstrate speedups ranging from 2.14× to 5.61×, reducing the maximum runtime from 42.12 minutes to 5.77 minutes across 45 test cases. These results highlight the scalability and computational efficiency of the proposed GPU–CPU concurrent solver framework.
\end{abstract}


\section{Introduction}

\subsection{Motivation}

Reliable and efficient operation of electric power systems is critical 
for modern infrastructure and economic activity \cite{b1}. The increasing scale, interconnection, and operational variability of contemporary grids have introduced significant challenges for system operators. These challenges require advanced decision-making tools capable of managing complex system interactions while ensuring secure and stable operation under a wide range of conditions. Recognizing these needs, the U.S. Department of Energy (DOE) introduced the Grid Optimization (GO) Challenges, led by the Advanced Research Projects Agency–Energy (ARPA-E) and Pacific Northwest National Laboratory (PNNL), to advance scalable optimization methods and computational strategies for large-scale power system applications \cite{b2}.

Within this context, optimal power flow (OPF) provides a mathematical framework for determining system operating points that satisfy physical and operational constraints while optimizing a specified objective. This process involves coordinating generator dispatch, load balancing, and energy storage operation to achieve efficient system performance \cite{b3, SCOPF-Javad}. To account for potential equipment failures and unforeseen events, security-constrained optimal power flow (SCOPF) extends the OPF formulation by incorporating contingency conditions directly into the optimization problem, thereby improving system preparedness and operational security \cite{b4}.

Historically, computational limitations restricted the use of detailed nonlinear models in SCOPF applications. As a result, simplified formulations such as DC optimal power flow (DCOPF) became widely adopted due to their reduced computational burden, despite neglecting important physical characteristics including voltage magnitude variations and reactive power flows \cite{b5}. Although advances in computing hardware and optimization algorithms have significantly expanded the feasibility of solving large-scale problems, practical solution methods for full AC-based formulations remain computationally demanding \cite{b2}.

\subsection{Existing Work on the GO Challenge} 
The GO Challenge focuses on solving large-scale unit commitment problems that incorporate nonlinear AC power flow constraints, making them particularly difficult to solve within practical time limits. These challenges have motivated the development of scalable optimization techniques capable of handling both discrete commitment decisions and nonlinear network constraints efficiently.

Several studies have investigated algorithmic enhancements to improve computational tractability. For instance, Holzer et al. \cite{refY} analyzed the performance of multiple optimization approaches submitted to GO Challenge 3, which involved multi-period unit commitment with AC network constraints and topology control. Their results demonstrated significant variability in solver performance depending on algorithm design and implementation, highlighting the importance of scalable and robust solution strategies.

In another line of work, Bienstock and Villagra \cite{refX} developed a cutting-plane-based linear relaxation for the AC optimal power flow (ACOPF) problem. Their approach generates valid outer-envelope cuts and applies refined cut selection techniques to strengthen relaxations and improve lower bound estimation, which is essential for addressing large-scale nonconvex optimization problems.

Research has also explored parallel and machine learning–inspired optimization methods. Chevalier \cite{ref4} proposed a parallel Adam-based optimization framework that accelerates convergence by distributing computations across multiple processing units, improving efficiency when solving high-dimensional optimization problems. Additionally, prior work in \cite{ref3} demonstrated that incorporating AC network constraints, reserve requirements, and unit commitment decisions within a unified formulation can produce reliable and feasible solutions for AC-based unit commitment problems.

In our earlier work in \cite{TPEC, IEEE-TIA, NAPS2025}, we introduced a decomposition-based solution approach tailored for the GO Challenge. The framework separates the problem into two stages: a DC-based formulation that determines binary unit commitment variables and ensures feasibility, followed by an AC-based formulation that resolves continuous variables using a nonlinear optimization solver. This sequential strategy improves computational efficiency while maintaining solution feasibility and accuracy for large-scale power system optimization problems.

\subsection{Contributions}

Building on our prior work \cite{TPEC26}, which demonstrated the use of GPUs to accelerate the unit commitment (UC) solving process, this study addresses the observation that not all scenarios benefit from GPU-based acceleration. In particular, certain problem instances remain more efficiently solved using CPU-based optimization algorithms due to differences in convergence behavior and numerical structure.

To address this limitation, this paper introduces a concurrent GPU--CPU solver framework that simultaneously executes GPU- and CPU-based optimization algorithms and terminates the slower solver once the faster solver converges. This approach systematically harnesses the complementary strengths of both architectures, improving robustness and reducing overall runtime across diverse UC instances.

In addition, the proposed framework is further strengthened through systematic tuning and refinement of the presolve strategy. The presolve procedure reduces problem size, improves numerical conditioning, and accelerates convergence for both solver pipelines.

The effectiveness of the proposed concurrent framework with tuned presolve is validated on a large-scale 6,049-bus system using the entire U.S. Department of Energy (DOE) power grid formulation. The results demonstrate significant speed improvements on all test cases without compromising solution quality.



\section{UC Solver Setup}

\subsection{Problem Formulation}
This study adopts the complete AC formulation defined by the DOE ARPA-E GO Challenge 3. The problem specification is extensive, including a 53-page primary document \cite{GO3} and 22 pages of supplementary material. It defines a large-scale optimization problem with detailed operational constraints, including real and reactive power limits at each bus, voltage bounds, zonal reserve requirements, unit commitment decisions, switching constraints, minimum uptime and downtime limits, startup restrictions, ramp rate limits, and energy production bounds over time. All constraints are enforced under both base-case and contingency conditions to ensure reliable system operation during component outages and unexpected disturbances.

\subsection{Baseline Model}

The baseline model is adopted from our previous work \cite{TPEC, IEEE-TIA, NAPS2025}. The model was tested on 300 different scenarios and achieved a scaled score of 0.98, indicating how close the solutions are to the best-known values. These scenarios were selected by domain experts, represent diverse and practical settings, and can be found in \cite{go_data}. Our recent work in \cite{TPEC26} achieved a significant speedup in computation time using GPU acceleration and the Primal-Dual Hybrid Gradient (PDHG) algorithm.
\subsection{Concurrent Solver Setup}
In this study, we propose a concurrent approach in which two solvers are run simultaneously on a single machine: when one solver converges, the other is terminated. One solver employs the PDHG algorithm to solve the relaxed linear problem using a GPU, while the second solver applies the barrier method to the same relaxed problem on CPU threads. Two solvers are used to leverage the advantages of both approaches, as in some cases the barrier method on CPU is faster than the GPU. 

Of the 56 available CPU cores, 48 cores are evenly allocated between the two solvers to enable fully concurrent execution and prevent resource contention (see Algorithm~\ref{alg:concurrent-solvers}). Since neither of the algorithms, PDHG nor barrier, returns a basic solution, both are followed by a crossover procedure and additional simplex iterations to refine the solution if necessary. This procedure provides an incumbent lower bound. Finally, the branch-and-bound method is applied to obtain the optimal solution within a specified optimality gap.

\begin{algorithm}[b]
\caption{Concurrent CPU--GPU Solver Framework}
\label{alg:concurrent-solvers}
\begin{algorithmic}[1]

\State Partition 48 CPU cores evenly between two solver instances.

\State Configure solver pipelines:
\Statex \hspace{1em}\textbf{Solver 1:}
\Statex \hspace{2em}Presolve $\rightarrow$ PDHG (GPU) $\rightarrow$ Crossover$\rightarrow$
\Statex \hspace{2em}Simplex (if needed)$\rightarrow$ Branch-and-Bound
\Statex \hspace{1em}\textbf{Solver 2:}
\Statex \hspace{2em}Presolve $\rightarrow$ Barrier $\rightarrow$ Crossover$\rightarrow$
\Statex \hspace{2em}Simplex (if needed)$\rightarrow$ Branch-and-Bound

\State Execute both solvers concurrently on the same optimization problem with an identical optimality gap tolerance.

\While{neither solver has converged}
    \State Wait 0.5 seconds.
    \If{either solver converges}
        \State Terminate the other solver.
        \State Return the converged solution.
    \EndIf
\EndWhile

\end{algorithmic}
\end{algorithm}

\subsection{Enhancing Solver Performance via Presolve}
\begin{table}[b]
\caption{Problem scale for different scheduling divisions (D1--D3), illustrating the large-scale nature of the optimization problem. Division~1 corresponds to a one-day horizon with 18 time steps, Division~2 to a two-day horizon with 48 time steps, and Division~3 to a one-week schedule with 42 time steps.}
\vspace{-.25cm}
\begin{center}
\begin{tabular}{|l|c|c|c|}
\hline
\textbf{Parameter} & \multicolumn{3}{c|}{\textbf{Divisions}} \\
\hline
 & \textbf{D1} & \textbf{D2} & \textbf{D3} \\
\hline
\hline
\# Continuous Variables & 1,741,233 & 4,111,861 & 3,599,256 \\
\hline
\# Binary Variables & 116,154 & 67,932 & 104,868 \\
\hline
\# Constraints & 1,580,545 &  3,826,499 & 3,348,659 \\
\hline
\end{tabular}
\label{tab01}
\end{center}
\end{table}

\begin{table}[h]
\centering
\caption{System configuration for experiments}
\label{tab:system_config}
\begin{tabular}{|l|l|}
\hline
\textbf{Hardware/Software} & \textbf{Details} \\ \hline
\hline
GPU & NVIDIA A800, 40 GB VRAM \\ \hline
CPU & Intel\textregistered{} Xeon\textregistered{} W9-3495X, 56 cores / \\
& 112 threads, 1.90–4.80 GHz \\ \hline
Memory (RAM) & 256 GB \\ \hline
OS & Linux 64-bit (via WSL on Windows) \\ \hline
Python & Version 3.11 \\ \hline
Gurobi & Version 13.0.1 beta\\ \hline
CUDA & Version 12.9 \\ \hline
Threads in Gurobi & 24 (default 32) \\ \hline
\end{tabular}
\end{table}
The PDHG algorithm can be tuned to reduce the time spent in the PDHG stage by, for example, setting the iteration limit to 100{,}000, the absolute feasibility tolerance to $10^{-5}$, and the convergence tolerance to $10^{-4}$. While these settings may produce a low-quality initial solution for the crossover procedure, which can increase the crossover solution time, they still reduce the overall solver runtime in most of the 45 scenarios tested in this study.

In this study, we adopt a more effective approach to reducing computation time by focusing on decreasing problem size and complexity. Specifically, we propose a more aggressive presolve strategy.

The presolve stage reduces the number of variables and constraints before the solver is called, which can improve the convergence of the optimization process. It is worth noting, however, that in some cases, smaller optimization problems can be more difficult to solve than larger ones due to the problem's structure and numerical stability issues, which may affect the convergence trajectory.

In this study, the number of presolve passes is doubled to explore additional opportunities to reduce problem complexity by enabling more simplifications. In this proposed setting, the presolve phase is also allowed to form the dual problem, enabling additional dual reductions that further decrease problem complexity. Although this increases the time spent in the presolve stage, the crossover procedure, simplex iterations, and branch-and-bound process are more likely to converge faster because they operate on a reduced problem size problem size.

\subsection{Test Instances}
Based on our study in \cite{TPEC26}, the most challenging system to solve is the 6049-bus network. In this study, we adopt this 6049-bus system to test the concurrent approach and the aggressive presolve technique. The data for Trials 2 and 3 are available from the ARPA-E dataset \cite{go_data}. In total, we use 45 test cases: 18 from Trial 3 and 27 from Trial 2.

Table~\ref{tab01} summarizes the problem scale for the different scheduling divisions (D1--D3). Each division represents a distinct temporal horizon: Division~1 corresponds to a one-day horizon with 18 time steps, Division~2 to a two-day horizon with 48 time steps, and Division~3 to a one-week schedule with 42 time steps. As shown in Table~\ref{tab01}, the optimization problems involve millions of continuous and binary variables and several million constraints, reflecting the large-scale and high-dimensional nature of these power system scheduling problems.

The size and complexity of these problems necessitate advanced computational strategies. 
By combining the concurrent approach, aggressive presolve, solver tuning, and hardware acceleration, we aim to improve convergence rates and reduce the overall computation time across all 45 test cases.

\subsection{Hardware}
The experiments were conducted on a Linux64 environment via Windows Subsystem for Linux (WSL) using an NVIDIA A800 GPU (40GB). Computations employed Gurobi 13.0.1 beta on a node with an Intel\textregistered{} Xeon\textregistered{} W9-3495X processor (56 cores, 112 threads, 1.90–4.80GHz) and 256GB of RAM. Each Gurobi solver instance was configured to use 24 threads, whereas the default setting is 32 threads.

GPU execution of the PDHG algorithm with Gurobi requires NVIDIA hardware (ideally an H100 GPU), CUDA 12.9, and a 64-bit Linux OS \cite{gurobi}. Although Gurobi supports Windows, GPU acceleration is currently limited to Linux systems. The complete hardware and software specifications used in this study are summarized in Table~\ref{tab:system_config}.

\section{Results}
\subsection{Impact of Aggressive Presolve on Solver Runtime}

\begin{figure}[ht]
    \centering
    \begin{subfigure}[t]{1\linewidth}
        \centering
        \includegraphics[width=1\linewidth]{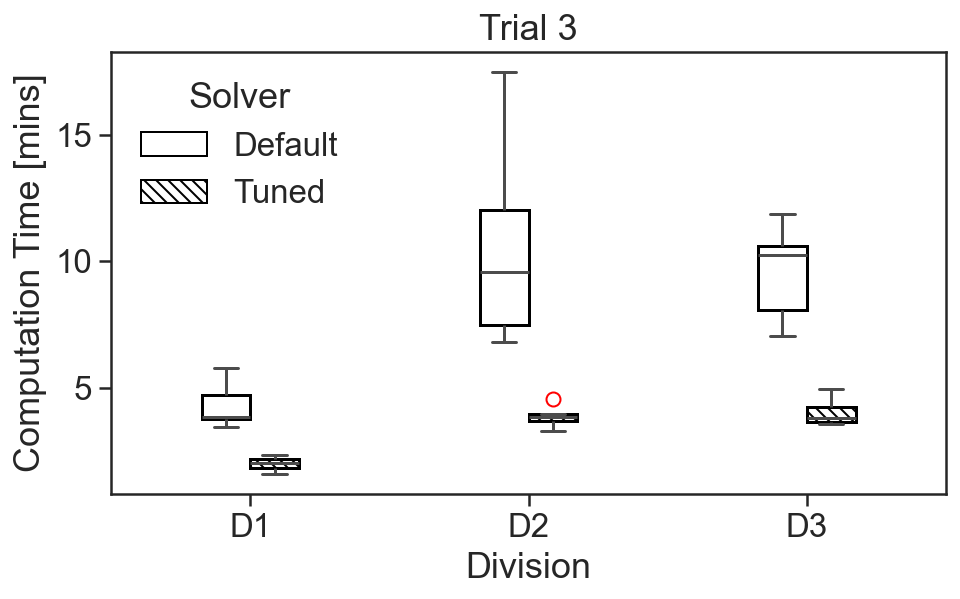}
        \caption{Trial 3}
        \label{fig:6079_trail3}
    \end{subfigure}

    \vspace{0.5cm}

    \begin{subfigure}[t]{1\linewidth}
        \centering
        \includegraphics[width=\linewidth]{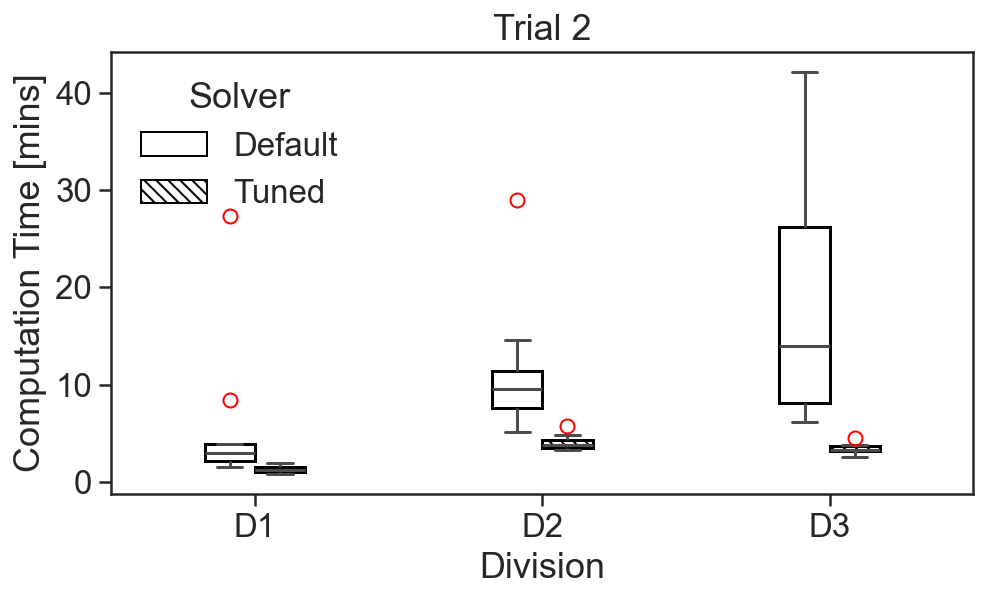}
        \caption{Trial 2}
        \label{fig:6079_trail2}
    \end{subfigure}

    \caption{Computation time distribution for the 6049-bus network. The tuned solver is consistently faster across scenarios and exhibits lower runtime variance for both Trail~3 and Trail~2 and for all divisions (D1-D3). Outliers are visualized using empty red circular markers.}
    \label{fig:6079_runtime}
\end{figure}

The results in Tables~\ref{tab:max_times_combined} and \ref{tab:avg_speedup_combined} demonstrate the significant impact of aggressive presolve on solver performance. Across all divisions (D1-D3) and Trials 2 and 3, the tuned solver (with enhanced presolve) consistently outperforms the untuned configuration. For example, in Trail~2, Division~3, the maximum runtime decreased from 42.12 minutes to 4.53 minutes, and the average runtime was reduced by 5.61× (Tables~\ref{tab:max_times_combined} and \ref{tab:avg_speedup_combined}). Similar improvements are observed across other divisions, with runtime speed-ups ranging from roughly 2× to 5.6×. The runtime distributions across all divisions and trials are shown using boxplots in Figure~\ref{fig:6079_runtime}. The results show that the tuned solver is consistently faster across scenarios and exhibits lower runtime variance for both Trials 2 and 3 across all divisions (D1–D3).

These results highlight that the aggressive presolve strategy effectively reduces problem size and complexity before solver execution. By eliminating redundant variables and constraints, enabling additional simplifications, and forming the dual problem for further reductions, the solver operates on a significantly smaller and more tractable problem. Consequently, not only do the PDHG, barrier, and crossover stages converge faster, but the simplex and branch-and-bound procedures also benefit from the reduced computational burden.

Overall, the combination of aggressive presolve, solver tuning, and concurrent computational techniques, particularly for large-scale instances such as the 6049-bus network, substantially enhances solver efficiency. This demonstrates that problem-specific presolve strategies are a critical component for tackling high-dimensional optimization problems in power system scheduling.

\begin{table}
\caption{\textcolor{black}{Maximum run-times for untuned (default) and tuned solvers for Trials 3 and 2. Arrows indicate whether tuned solver is faster ($\downarrow$, blue) or slower ($\uparrow$, red) than untuned. Time in minutes.}}
\centering
\setlength{\tabcolsep}{6pt}
\begin{tabular}{|l|c|c|c|}
\hline
\textbf{Trail} & \textbf{Division (D)} & \multicolumn{2}{c|}{\textbf{Maximum Runtime [min]}} \\
\cline{3-4}
 &  & \textbf{Default} & \textbf{Tuned} \\
\hline
\hline
Trail 3 & D1 & 5.78 & 2.34 \textcolor{blue}{$\downarrow$} \\
        & D2 & 17.45 & 4.54 \textcolor{blue}{$\downarrow$} \\
        & D3 & 11.87 & 4.96 \textcolor{blue}{$\downarrow$} \\
\hline
Trail 2 & D1 & 27.35 & 1.91 \textcolor{blue}{$\downarrow$} \\
        & D2 & 29.00 & 5.77 \textcolor{blue}{$\downarrow$} \\
        & D3 & 42.12 & 4.53 \textcolor{blue}{$\downarrow$} \\
\hline
\end{tabular}
\label{tab:max_times_combined}
\end{table}

\begin{table}
\caption{\textcolor{black}{Average runtime for untuned (default) and tuned solvers and resulting speed-up for Trials 3 and 2. Colors indicate whether tuning improves performance (blue) or degrades it (red). Time in minutes.}}
\centering
\setlength{\tabcolsep}{6pt}
\begin{tabular}{|l|c|c|c|c|}
\hline
\textbf{Trail} & \textbf{Division (D)} & \multicolumn{2}{c|}{\textbf{Avg. Runtime [min]}}  & \textbf{Avg.} \\
\cline{3-4}
 &  & \textbf{Default} & \textbf{Tuned} & \textbf{Speed-up} \\
\hline
\hline
Trail 3 & D1 & 4.28 & 2.01 & \textcolor{blue}{$2.14\texttt{x}$} \\
        & D2 & 10.51 & 3.85 & \textcolor{blue}{$2.73\texttt{x}$} \\
        & D3 & 9.58 & 4.02 & \textcolor{blue}{$2.38\texttt{x}$} \\
\hline
Trail 2 & D1 & 5.93 & 1.27 & \textcolor{blue}{$4.68\texttt{x}$} \\
        & D2 & 11.45 & 4.05 & \textcolor{blue}{$2.82\texttt{x}$} \\
        & D3 & 18.98 & 3.38 & \textcolor{blue}{$5.61\texttt{x}$} \\
\hline
\end{tabular}
\label{tab:avg_speedup_combined}
\end{table}

\subsection{Assessment of Solution Quality}

Table~\ref{tab:solution_quality_distribution} reports the solution quality score, defined as the ratio between the objective value obtained by the tuned solver and that obtained by the default solver. Since the objective function represents market surplus, which is maximized, a ratio greater than one indicates an improvement in solution quality.

Out of the 45 scenarios considered from Trails 2 and 3, the proposed tuned solver maintains the solution quality in 30 scenarios. In the remaining 15 scenarios, there is a very small reduction in solution quality; however, the score remains above 0.9999999 in all cases. The average score across the 45 scenarios is 0.99999998, indicating that the proposed tuned solver effectively preserves the solution quality.

This demonstrates that solver tuning can achieve substantial speed improvements without compromising solution optimality or feasibility.

\begin{table}[b]
\caption{Distribution of solution quality for Tuned vs default solvers. Score $=$ Tuned $/$ Untuned objective. Out of 45 scenarios from Trails~2 and~3, 30 scenarios achieve a ratio of at least~1, while 15 scenarios have ratios slightly below~1 but above~0.9999999, indicating that the tuned solver preserves solution quality across all tested cases.}
\centering
\begin{tabular}{|l|c|}
\hline
\textbf{Score Range} & \textbf{Number of Scenarios} \\
\hline
\hline
 Score $<0.9999999$& \textcolor{blue}{0} \\
\hline
$0.9999999 \leq$ Score $< 1$ & \textcolor{red}{15} \\
\hline
\textbf{Score $\geq 1$} & \textcolor{blue}{30} \\
\hline
\hline
\textbf{Average Score} & $\approx$ \textcolor{red}{0.99999998} \\
\hline
\end{tabular}
\label{tab:solution_quality_distribution}
\end{table}

\subsection{Impact of Solver Tuning on Runtime Breakdown}

\begin{figure*}[ht]
    \centering
    \begin{subfigure}[b]{0.48\linewidth}
        \centering
        \includegraphics[width=1\linewidth]{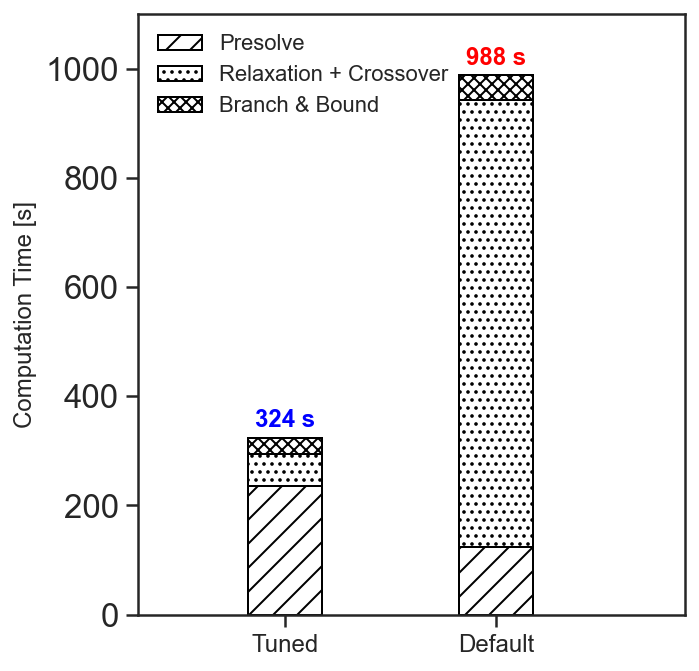}
        \caption{Scenario 31 from D2}
        \label{fig:Time_breakdown1}
    \end{subfigure}
    \hfill
    \begin{subfigure}[b]{0.48\linewidth}
        \centering
        \includegraphics[width=1\linewidth]{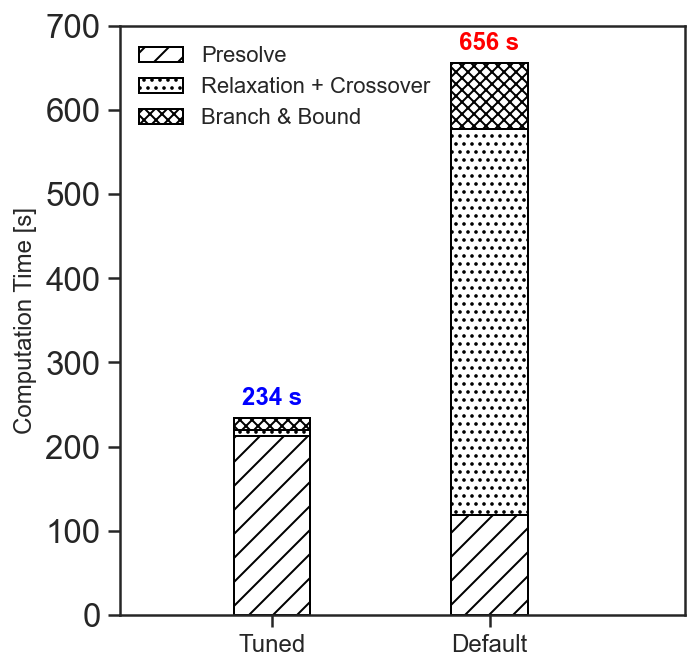}
        \caption{Scenario 42 from D3}
        \label{fig:Time_breakdown2}
    \end{subfigure}
    \caption{Comparison of computation times between tuned and default solvers across two random scenarios. Each bar represents the total runtime, decomposed into three stages: (1) presolve, (2) relaxation and crossover, and (3) branch-and-bound. The tuned solver reduces model complexity by spending more time in the presolve stage. Although the time spent in presolve nearly doubles, the overall computation time decreases significantly due to faster convergence in the relaxation, crossover, and branch-and-bound stages.}
    \label{fig:Time_breakdown_combined}
\end{figure*}

In Scenario 31 from D2 (Fig.~\ref{fig:Time_breakdown1}), the presolve stage of the tuned solver consumes 235s, compared to 123s for the untuned solver. Although the presolve time nearly doubles, the relaxation and crossover stage is significantly reduced, requiring only 60s for the tuned solver versus 820s for the untuned solver. Similarly, the branch-and-bound stage decreases from 45s to 29s. As a result, the total runtime of the tuned solver is 324s, compared to 988s for the untuned solver, corresponding to an overall speedup of approximately 3×.

In Scenario 42 from D3 (Fig.~\ref{fig:Time_breakdown2}), the tuned solver again spends more time in the presolve stage (213s) compared to the untuned solver (119s). This additional presolve effort leads to substantial reductions in the subsequent stages: the relaxation and crossover time decreases from 459s to 7s, and the branch-and-bound time is reduced from 78s to 14s. Consequently, the total runtime drops from 656s for the untuned solver to 234s for the tuned solver, yielding a speedup of approximately 2.8×.

For the CPU-based solver, the relaxation and crossover stage consists of three components: solving the continuous relaxation using a barrier method, performing the crossover procedure, and, when necessary, executing additional simplex iterations. In some cases, the crossover step terminates with elevated dual infeasibility, which requires simplex iterations to recover a high-quality basic feasible solution. Accordingly, the reported relaxation and crossover time includes the barrier solve, the crossover process, and any subsequent simplex iterations used to clean the solution. The bounds obtained from this stage are then used to guide the branch-and-bound process in determining optimality.

In contrast, for the GPU-based solver, the continuous relaxation is solved using the PDHG algorithm, which replaces the barrier method. PDHG is particularly well suited for large-scale problems and typically achieves faster convergence, leading to a significant reduction in the relaxation time and, consequently, the overall computation time.

\subsection{GPU–CPU Concurrency Under Solver Tuning}
\begin{table}[b]
\centering
\small
\caption{Number of scenarios (out of 45) solved first by each solver in concurrent CPU–GPU mode under different solver settings.}
\begin{tabular}{|l|l|l|}
\hline
\textbf{Settings} & \textbf{GPU Solver First} & \textbf{CPU Solver First} \\ 
\hline
\hline
Default & 33 \textcolor{blue}{(73\%)} & 12 (27\%) \\ 
Tuned & 0 (0\%) & 45 \textcolor{blue}{(100\%)} \\ 
\hline
\end{tabular}
\label{T2_sync}
\end{table}
To recall, in the proposed setup, the GPU-based and CPU-based solvers are executed concurrently, and the run is terminated as soon as one solver converges. 

Without solver tuning, the GPU-based solver is faster than the CPU-based solver in most cases. As shown in Table~\ref{T2_sync}, the GPU-based solver finishes first in 33 of the 45 test cases, whereas the CPU-based solver finishes first in only 12. This result suggests that the GPU architecture provides a clear computational advantage for the original large-scale problem due to its ability to exploit a high degree of parallelism.
However, under solver tuning, the CPU-based solver consistently converges earlier than the GPU-based solver. This behavior is primarily due to the aggressive presolve performed by the tuned solver, which substantially reduces the problem size and model complexity. Consequently, the remaining optimization problem becomes sufficiently compact that the GPU-based relaxation solver no longer offers a computational advantage over the CPU-based solver, thereby diminishing the effectiveness of GPU acceleration in this setting. These observations highlight that solver tuning can outweigh hardware-based acceleration for certain large-scale optimization problems.
\subsection{Impact of Solver Tuning on Model Size}
For large-scale problems (D2 and D3), the tuned solver aggressively reduces both constraints and variables, eliminating nearly all binary variables and collapsing the problem to a small continuous core (Table~\ref{tab:presolve_comparison}). 
Specifically, for D2, the tuned solver reduces constraints from 3,826,499 to 38,238 ($\approx 99\%$) and variables from 4,115,635 to 44,822 ($\approx 98.9\%$). Similarly, for D3, constraints drop from 3,348,659 to 38,238 ($\approx 98.9\%$) and variables from 3,603,030 to 44,822 ($\approx 98.8\%$). This substantial reduction dramatically accelerates subsequent stages, including relaxation, crossover, simplex, and branch-and-bound.
\begin{table*}[h!]
\caption{Total number of variables (continuous and binary) and constraints before and after presolve for default and tuned solvers across Divisions D1--D3. Percentages indicate reduction relative to the original model size.}
\centering
\setlength{\tabcolsep}{4pt}
\begin{tabular}{|l|l|l|l|l|l|l|l|}
\hline
\textbf{Division} & \textbf{Solver} 
& \multicolumn{3}{c|}{\textbf{Constraints}} 
& \multicolumn{3}{c|}{\textbf{Variables (Total)}} \\
\cline{3-8}
 & \textbf{Settings} & Before & After & \% Reduction & Before & After & \% Reduction \\
\hline
\hline
D1 & Default & 1,580,545 & 679,208 & 57.0\% & 1,857,387 & 791,966 & 57.4\% \\
D1 & Tuned   & 1,580,545 & 679,208 & 57.0\% & 1,857,387 & 795,966 & 57.1\% \\
\hline
D2 & Default & 3,826,499 & 1,108,844 & 71.0\% & 4,115,635 & 1,299,780 & 68.4\% \\
D2 & Tuned   & 3,826,499 & 38,238 & \textcolor{blue}{99.0\%} & 4,115,635 & 44,822 & \textcolor{blue}{98.9\%} \\
\hline
D3 & Default & 3,348,659 & 1,032,372 & 69.2\% & 3,603,030 & 1,210,140 & 66.4\% \\
D3 & Tuned   & 3,348,659 & 38,238 & \textcolor{blue}{98.9\%} & 3,603,030 & 44,822 & \textcolor{blue}{98.8\%} \\
\hline
\end{tabular}
\label{tab:presolve_comparison}
\end{table*}

For the smaller D1 problem, the tuned presolve retains slightly more variables than the untuned solver, reducing variables from 1,857,387 to 795,966 ($\approx 57.1\%$) compared to 791,966 ($\approx 57.4\%$) for the untuned solver. This behavior is due to additional presolve passes and dual reductions that preserve certain variables and slack structures, improving convergence stability and solution quality. Consequently, while the absolute reduction is smaller for D1, the solver ensures faster and high-quality solutions.

\subsection{Impact of Threads on Solver Speed}
Although the machine used for the experiments has 112 threads, only 48 threads were actively used and evenly partitioned between the two solver instances. Increasing the number of threads beyond this level did not noticeably improve convergence speed and, in some cases, even worsened computation time. This behavior is likely due to memory bandwidth limitations and increased thread synchronization overhead, which can reduce the efficiency of parallel execution when too many threads are used. For larger power grid instances, however, the optimal number of threads may differ, as a higher computation-to-communication ratio could make additional threads more beneficial. 
\vspace{-.1cm}
\section{Conclusions}
This paper developed a GPU–CPU concurrent solver for large-scale unit commitment problems. By running different solvers concurrently to leverage the strengths of both architectures and by tuning the presolve phase, the framework achieved significant runtime reductions, with the maximum runtime decreasing from 42.12 minutes to 5.77 minutes and speedups ranging from 2.14× to 5.61× across 45 test cases on a 6,049-bus system. These results demonstrate that the proposed approach efficiently handles very large optimization problems, providing scalable and practical solution strategies for real-time power system operations while maintaining the quality and feasibility of the computed solutions. Future work will explore extending both the concurrent solvers strategy and the proposed aggressive presolve approach to additional classes of large-scale power system optimization problems. 
\vspace{-.1cm}
\section*{Acknowledgment}
The authors acknowledge the use of artificial intelligence tools to improve the readability of this manuscript. All AI-assisted edits were verified and revised by the authors.
The authors also thank the support provided by the U.S. ARPA-E (Grant \#DE-AR0001646) and the NSF (Grants \#2442420 and \#2313768).

\end{document}